\documentclass[12 pt]{article}
\usepackage{amssymb}
\usepackage{amssymb,amsbsy,float}
\usepackage{graphicx}
\usepackage{amsfonts}
\usepackage{latexsym}
\usepackage{bibtopic}

\newtheorem{theorem}{Theorem}

\newtheorem{corollary}[theorem]{Corollary}

\newtheorem{proposition}[theorem]{Proposition}

\newtheorem{lemma}[theorem]{Lemma}

\newenvironment{proof}[1][Proof]{\textbf{#1.} }{\ \rule{0.5em}{0.5em} \vspace{0.1 cm}}
\newcommand{\bx}{\mbox{$\textbf{x}$}}

\newcommand{\R}{\mbox{$\mathbb{R}$}}%
\newcommand{\N}{\mbox{$\mathbb{N}$}}%

\begin{document}
\title{Equilibrium payoffs in finite games}
\author{Ehud Lehrer, Eilon Solan and Yannick Viossat\thanks{E. Lehrer and E. Solan: School of Mathematical Sciences, Tel-Aviv University, Israel; Y. Viossat: CEREMADE, Universit\'e Paris-Dauphine, France. E-mail to: viossat@ceremade.dauphine.fr.} \thanks{Y. Viossat thanks participants to the game theory seminar of the Institut Henri Poncar\'e, Paris, and to the ``Communication and networks in games" workshop in Valencia, Spain, for which PICASSO funding is gratefully acknowledged.}}
\maketitle
\begin{abstract}
We study the structure of the set of equilibrium payoffs in finite games, both for Nash equilibrium and correlated equilibrium. A nonempty subset of $\R^2$ is shown to be the set of Nash equilibrium payoffs of a bimatrix game if and only if it is a finite union of rectangles. Furthermore, we show that for any nonempty finite union of rectangles $U$ and any polytope $P \subset \R^2$ containing $U$, there exists a bimatrix game with $U$ as set of Nash equilibrium payoffs and $P$ as set of correlated equilibrium payoffs. The $n$-player case and the robustness of this result to perturbation of the payoff matrices are also studied.\\


\end{abstract}
\section{Introduction}
It is well known that the set of correlated equilibrium payoffs of a $n$-player finite game is a polytope 
in $\R^n$ (Aumann, 1974). We show that, conversely, any polytope in $\R^n$ is the set of correlated 
equilibrium payoffs of a $n$-player game. Moreover, for any $n$-player game $G$ and any polytope 
$P$ in $\R^n$ containing the correlated equilibrium payoffs of $G$, there exists a $n$-player game 
with the same set of Nash equilibrium payoffs as $G$ and $P$ as set of correlated equilibrium payoffs. 
If $G$ is a two-player game, it suffices to assume that $P$ contains the Nash equilibrium payoffs of $G$ 
(but not necessarily its correlated equilibrium payoffs). 
We also show that a subset of $\R^2$ is the set of Nash equilibrium payoffs of a finite game if and only 
if it is a finite union of rectangles. Together with the previous result, this implies that: for any nonempty 
finite union of rectangles $U$ and any polytope $P$ containing $U$, there exists a bimatrix game 
with $U$ as set of Nash equilibrium payoffs and $P$ as set of correlated equilibrium payoffs.

Finally, we study the robustness of our results to perturbation of the payoff matrices. Since almost all games have a finite set of equilibria, all we can hope to show is that, for any finite set $U$ in $\R^n$ and any polytope $P$ in $\R^n$ containing $U$, there exists an open set of games whose set of Nash equilibrium payoffs is arbitrarily close to $U$ and 
whose set of correlated equilibrium payoffs is arbitrarily close to $P$. We show that this is indeed the case.\\ 

\emph{Notations and definitions}:  throughout, $n$ denotes an integer weakly greater than $2$. Let $G$ be a $n$-player game with pure strategy set $S_i$ and payoff function $u_i : S:=\times_{1 \leq i \leq n} S_i \to \R$ for player $i$. Let $S_{-i}:=\times_{j \neq i} S_j$. A correlated equilibrium of $G$ (Aumann, 1974) is a probability distribution $\mu$ on the set of strategy profiles $S$ of $G$ such that, for every player $i$ in $\{1,..,n\}$ and every pure strategy $s_i$ in $S_i$:
\begin{equation}
\label{eq:defcor}
\forall t_i \in S_i, \sum_{s_{-i} \in S_{-i}} \mu(s_i,s_{-i}) [u_i(s_i,s_{-i})-u_i(t_i,s_{-i})] \geq 0
\end{equation}
The set of correlated equilibria of $G$ is a polytope, which contains the Nash equilibria. An extreme correlated equilibrium is an extreme point of this polytope. 

Let $u_i(\mu):=\sum_{s \in S} \mu(s) u_i(s)$ denote the average payoff of player $i$ in the correlated equilibrium $\mu$. The set of correlated equilibrium payoffs of $G$ is the set of $n$-tuples $(u_1(\mu),...,u_n(\mu))$ where $\mu$ is a correlated equilibrium of $G$. This is a polytope in $\R^n$. It is denoted by $CEP(G)$. Similarly, $NEP(G)$ denotes the set of Nash equilibrium payoffs of $G$, and we may write NE and CE for ``Nash equilibrium" and ``correlated equilibrium", respectively. Finally, if $A$ is a subset of $\R^n$ then $Conv(A)$ denotes its convex hull.

\section{Main results}
\label{sec:main}
We begin with a result on correlated equilibria: 
\begin{proposition}
\label{prop:cor}
For any polytope $P$ in $\R^n$, there exists a
$n$-player game whose set of correlated equilibrium payoffs (resp.
convex hull of the set of Nash equilibrium payoffs) is $P$.
\end{proposition}
\begin{proof}
We prove the result for two-player games; see the appendix for the $n$-player case. Let $P$ be a polytope in
$\R^2$. Let $(x_1,y_1)$,...,$(x_m, y_m)$ be its extreme points.
Assume that for every $i$ in $\{1,2,..,m\}$, $x_i$ and $y_i$ are
strictly positive. This is without loss of generality, since adding
a constant to all payoffs of a game does not change the set of Nash
equilibria and correlated equilibria. Consider the $m \times m$ game
with payoff matrix $(a_{ij},b_{ij})_{1 \leq i,j \leq m}$ such that,
for every $i$ in $\{1,2,...,m\}$, $a_{ii}=a_{mi}=x_i$,
$b_{ii}=b_{im}=y_i$, and all other payoffs are zero. For $m=4$ this
gives :
\begin{equation}
\label{game:cor}
\left(
  \begin{array}{llll}
    x_1,y_1 & 0,0       &  0,0    & 0,y_1\\
    0,0       & x_2,y_2 &   0,0    & 0,y_2 \\
    0,0       & 0,0       & x_3,y_3& 0,y_3\\
    x_1,0     & x_2,0     &  x_3,0   & x_4,y_4
  \end{array}
\right)
\end{equation}
Clearly, any diagonal square corresponds to a pure Nash equilibrium.
Furthermore, choosing the last row (column) is a weakly dominant
strategy for the row (column) player, and it gives a strictly higher
payoff than choosing row (column) $i \neq m$ whenever the column
(row) player does not choose column (row) $i$. It follows that in
every correlated equilibrium, the probability of every off-diagonal
square is zero, hence there are no other Nash equilibria or extreme
correlated equilibria than the diagonal squares. Therefore, the
convex hull of the set of Nash equilibrium payoffs coincides with
the set of correlated equilibrium payoffs and is equal to $P$.
\end{proof}

To state our next result, we first need a definition:  let us say that a subset of $\R^2$ is a \emph{rectangle} if it is of the
form $[a,b] \times [c,d]$, for some real numbers $a$, $b$, $c$, $d$, with $a \leq b$, $c\leq d$.
\begin{proposition}
\label{prop:Nash}
In any bimatrix game, the set of Nash equilibrium payoffs is a
finite union of rectangles. Conversely, for any nonempty finite union of rectangles $U$, there exists a
bimatrix game whose set of Nash equilibrium payoffs is $U$. Moreover, for any bimatrix game $G$, there exists a bimatrix game $G'$ with $NEP(G')=NEP(G)$ and $CEP(G')=Conv(NEP(G))$. 
\end{proposition}
\begin{proof}
Consider a bimatrix game with pure strategy set $I$ for player $1$ and $J$ for player $2$. Let $NEP \subset \R^2$ denote its set of Nash equilibrium payoffs. For any subsets $I'$ of $I$ and $J'$ of $J$, let $NEP(I',J')$ denote the set of payoffs of Nash equilibria with support $I' \times J'$. For any subset $A$ of $\R^2$, let $\bar{A}$ denote the closure of $A$. We have 
$$NEP=\bigcup_{I' \subset I, J' \subset J} NEP(I',J')= \bigcup_{I'\subset I, J'\subset J} \overline{NEP(I',J')}$$
where the second equality holds because $NEP$ is closed. Thus, to show that $NEP$ is a finite union of rectangles, it suffices to show that $\overline{NEP(I',J')}$ is a rectangle. This follows easily from the fact that, in bimatrix games, Nash equilibria with the same support are exchangeable.

Conversely, let $m \in \N^*$ and for $1 \leq i \leq m$, let $a_i$, $b_i$, $c_i$, $d_i$ be real numbers. Let
$U= \bigcup_{1 \leq i \leq m} [a_i,b_i] \times [c_i,d_i]$. Assuming w.l.o.g. that the numbers $a_i$, $b_i$, $c_i$, $d_i$ are all positive, we build below a bimatrix game with $U$ as set of Nash equilibrium payoffs. Consider first the game with payoff matrices:
\[(A_i,B_i)= \left(
                  \begin{array}{cc}
                    a_i, c_i & b_i, c_i \\
                    a_i, d_i & b_i, d_i \\
                  \end{array}
                \right)
\]
In this game, a player does not influence its own payoffs and the set of Nash equilibrium payoffs is $[a_i,b_i] \times [c_i,d_i]$. Let
$$(A_i,0)=\left(
                  \begin{array}{cc}
                    a_i,0 & b_i,0 \\
                    a_i,0 & b_i,0 \\
                  \end{array}
                \right) \quad , \quad (0, B_i)=\left(
                  \begin{array}{cc}
                    0, c_i & 0, c_i \\
                    0, d_i & 0, d_i \\
                  \end{array}
                \right)$$
and consider the game built by blocks:
\begin{equation}\label{game:Nash}
\left(
  \begin{array}{ccccc}
    (A_1,B_1) & 0         & .. &                0 & (0, B_1) \\
    0         & (A_2,B_2) & .. &                0 & (0, B_2) \\
    ..        & ..        & .. & ..               & ..        \\
    0         & 0         & .. & (A_{m-1},B_{m-1})& (0, B_{m-1}) \\
    (A_1,0)   & (A_2,0)   & .. & (A_{m-1},0)      & (A_m,B_m) \\
  \end{array}
\right)
\end{equation}
where an isolated $0$ represents a $2 \times 2$ block of payoffs $(0,0)$. This game has the same structure as (\ref{game:cor}), but the payoffs $x_i$ and $y_i$ have been replaced by the blocks $A_i$ and $B_i$, respectively.   Any mixed strategy profile with support in one of the blocks $(A_i,B_i)$ is a Nash equilibrium. Furthermore, it is easy to prove along the lines of the proof of proposition \ref{prop:cor} that there are no other Nash equilibria and that the set of correlated equilibria is the convex hull of the set of Nash equilibria. It follows that: first, the set of NE payoffs of (\ref{game:Nash}) is equal to $U$, which shows that any finite union of rectangles is the set of NE payoffs of a bimatrix game; and second, the set of CE payoffs of (\ref{game:Nash}) is equal to $Conv(U)$. Since the set of NE payoffs of a bimatrix game is always a finite union of rectangles, this implies that for any bimatrix game $G$, there exists a bimatrix game $G'$ (of type (\ref{game:Nash})) such as $NEP(G')=NEP(G)$ and $CEP(G')=Conv(NEP(G))$. \end{proof}

\begin{lemma}
\label{lm:add-cor}
Let $(x_1,x_2,...,x_n) \in \R^n$. For any $n$-player game $G$, there exists a
$n$-player game with the same set of Nash equilibrium payoffs as $G$
and whose set of correlated equilibrium payoffs is the convex hull
of $(x_1,x_2,...,x_n)$ and of the set of correlated equilibrium payoffs of $G$.
\end{lemma}

\begin{proof}
We prove the result in the two-player case. For the $n$-player case, see the appendix. Let $G$ be a two-player $m_1 \times m_2$
game and $(x, y) \in \R^2$. Assume w.l.o.g. that
$x$, $y$, and all the payoffs in $G$ are positive. Consider
the $(3 + m_1) \times (3 + m_2)$ game
\begin{equation}
\label{game:lm}
\Gamma=\left(\begin{array}{c|c}
       \begin{array}{ccc}
         0,0        & x+1, y-1  & x-1, y+1 \\
         x-1, y+1   & 0,0       & x+1, y-1 \\
         x+1, y-1   & x-1, y+1  & 0,0
       \end{array}
                            & [0,y]
        \\
        \hline
       [x,0] & G
     \end{array}\right)
\end{equation}
where $[x,0]$ denotes a block of payoffs $(x,0)$
of appropriate size (same notations for player $2$). 
Let $\nu$ denote the correlated strategy putting
probability $1/6$ on every off-diagonal square of the top-left
block, and probability $0$ on every other square of the whole payoff
matrix. Clearly, $\nu$ is a correlated equilibrium of $\Gamma$, with
payoff $(x,y)$, and every correlated equilibrium of $G$ induces a
correlated equilibrium of $\Gamma$. We claim that any correlated
equilibrium of $\Gamma$ is a convex combination of $\nu$ and of a
correlated equilibrium of $G$. 

To see this, let $\mu$ be a correlated equilibrium of $\Gamma$. Clearly, $\mu_{12} \geq \mu_{13}$, 
otherwise player $1$ would have an incentive to deviate from his first to his last strategy 
(recall that all payoffs in $G$ are positive). Repeating this reasoning with other strategies 
and with player $2$ leads to the chain of inequalities
$$\mu_{12} \geq \mu_{13} \geq \mu_{23} \geq \mu_{21} \geq \mu_{31} \geq \mu_{32} \geq
\mu_{12}$$
Since the first and last terms are equal, this is a chain of equalities, hence $\mu$ puts equal weight 
on all off-diagonal squares of the top-left block. It is then easy to see that $\mu$ puts probability zero 
on the diagonal of the top-left block as well as on the top-right and bottom-left blocks. This implies that $\mu$ is a convex
combination of $\nu$ and of a correlated equilibrium of $G$, proving the claim. 

It follows that : (i) a fortiori, any CE payoff of $\Gamma$ is a convex combination of the payoff of $\nu$ and of a correlated equilibrium payoff of $G$; 
and (ii) $\Gamma$ and $G$ have the same Nash equilibria, hence a fortiori the same NE payoffs.  This concludes the proof.
\end{proof}

\begin{proposition}\label{prop:NashAndCor}
For any $n$-player game $G$ and for any polytope $P$ in $\R^n$ containing
$CEP(G)$, there exists a game with the same set of Nash equilibrium
payoffs as $G$ and $P$ as set of correlated equilibrium payoffs.
\end{proposition}
\begin{proof}
Let $P$ be a polytope containing $CEP(G)$, with $q$ extreme points. Applying iteratively lemma \ref{lm:add-cor} ($q$ times), we obtain a game with the same NE payoffs as $G$ and whose set of CE payoffs is the convex hull of $P$ and $CEP(G)$, hence is equal to $P$.\end{proof}

In the two-player case, we get a stronger result:
\begin{corollary}
\label{cor:2player}
For any bimatrix game $G$ and for any polytope $P$ in $\R^n$ containing
$NEP(G)$, there exists a game with the same set of Nash equilibrium
payoffs as $G$ and $P$ as set of correlated equilibrium payoffs.\footnote{The difference with the $n$-player case in that we only assume $P \supset NEP(G)$ instead of $P \supset CEP(G)$.} 
\end{corollary}
\begin{proof}
By proposition  \ref{prop:Nash}, there exists a game $G'$ with $NEP(G')=NEP(G)$ and $CEP(G')=Conv(NEP(G))$. Applying proposition \ref{prop:NashAndCor} to $G'$ gives the result.
\end{proof}

In light of proposition  \ref{prop:Nash}, corollary \ref{cor:2player} may also be stated as follows : for any nonempty finite union of rectangles $U$ in $\R^2$ and for any polytope $P$ in $\R^2$ containing $U$, there exists a two-player game with $U$ as set of Nash equilibrium payoffs and $P$ as set of correlated equilibrium payoffs. Proposition \ref{prop:Nash} and corollary \ref{cor:2player} fully characterize the possible shapes of the set of Nash equilibrium payoffs and correlated equilibrium payoffs for two-player games. For $n$-player games, the general shape of the set of Nash equilibrium payoffs remains to be understood. Furthermore, we do not know whether, for any $n$-player game $G$, there exists a $n$-player game with the same set of Nash equilibrium payoffs as $G$ and whose set of correlated equilibrium payoffs is the convex hull of its Nash equilibrium payoffs. 
\section{Genericity}
The games used above are highly non generic. For instance, a small perturbation of the payoffs of (\ref{game:cor}) is enough to eliminate all its Nash and correlated equilibria but one. This raises the issue of the robustness of our results. Ideally, to show that, for instance, corollary \ref{cor:2player} is robust, one would like to show that for any nonempty finite union of rectangles $U$ and for any polytope $P$ in $\R^2$ containing $U$, there exists an open set of games whose set of Nash equilibrium payoffs is ``close" to $U$ and whose set of correlated equilibrium payoffs is ``close" to $P$. This is hopeless however, since almost all games have a finite set of Nash equilibria. Thus, all we can hope to prove is the same result but with $U$ being a finite set. This is the object of this section.

We first need some definitions.  Let $\epsilon>0$. For any $\bx$ in $\R^n$, let $||\bx||=\max_{1\leq i \leq n} |x_i|$. 
Let $A$ and $A'$ be subsets of $\R^n$. Recall that $A$ and $A'$ are $\epsilon$-close in the Hausdorff distance sense if
\begin{equation}
\forall \bx \in A, \exists \bx' \in A', ||\bx-\bx'|| < \epsilon
\end{equation}
and
\begin{equation}
\forall \bx' \in A', \exists \bx \in A, ||\bx-\bx'|| < \epsilon
\end{equation}
Let $G$ be a finite game with payoff function $u_i$ for player $i$ and let $\alpha >0$. The open ball of center $G$ and radius $\alpha$, denoted by $B(G,\alpha)$, is the set of all games $G'$ with the same sets of players and strategies as in $G$ and such that for every player $i$ and every pure strategy profile $s$, $|u'_i(s)-u_i(s)| < \alpha$, where $u'_i$ is the payoff function of player $i$ in $G'$. A set of games $\Sigma$ is open if for every game $G$ in $\Sigma$, $\Sigma$ contains an open ball of center $G$ and positive radius.
\begin{proposition}
\label{prop:gen}
Let $U$ be a finite set in $\R^n$. Let $P \subset \R^n$ be a polytope containing $U$.
There exists a (nonempty) open set of $n$-player games whose set of Nash equilibrium payoffs is $\epsilon$-close to $U$ and
whose set of correlated equilibrium payoffs is $\epsilon$-close to $P$.
\end{proposition}
\begin{proof}
We prove the result for two-player games. The proof for $n$-player games is similar (see the appendix).
Let $U=\{(x_1,y_1),...,(x_m, y_m)\}$, let $P \subset \R^2$ be a polytope containing $U$, with vertices $(x'_1,y'_1)$,...,$(x'_q, y'_q)$. Assume w.l.o.g. 
that, for all $i$ in $\{1,..,q\}$,  $x_i$, $y_i$, $x'_i$ and $y'_i$ are all positive. Let $G_{\alpha}$ denote the $m \times m$ game
with payoff matrix $(a_{ij},b_{ij})_{1 \leq i,j \leq m}$ such that:
for every $i$ in $\{1,2,...,m\}$, $a_{ii}=x_i$ and $b_{ii}=y_i$; for every $i$ in $\{1,2,...,m-1\}$, $a_{mi}=x_i-\alpha$ and $b_{im}=y_i-\alpha$; and all other payoffs are zero. For $m=4$ this
gives:
$$
G_{\alpha} = \left(
  \begin{array}{cccc}
    x_1,y_1     & 0 , 0           &   0,0         & 0,y_1-\alpha\\
    0 , 0         & x_2,y_2       &   0,0         & 0,y_2-\alpha \\
    0,0         & 0,0           & x_3,y_3       & 0,y_3-\alpha\\
   x_1-\alpha,0 & x_2-\alpha,0  & x_3-\alpha,0  & x_4,y_4
  \end{array}
\right)
$$
Thus, $G_0$ is the game used in the proof of proposition \ref{prop:cor} and $NEP(G_0)=U$. 
For every $(x,y)$ in $\R^2$, let $C(x,y)$ denote the game corresponding to the top-left block of (\ref{game:lm}).
\begin{equation}
\label{eq:Cxy}
C(x,y)=\left(\begin{array}{ccc}
         0,0           & x+1, y-1 & x-1, y+1\\
         x-1, y+1  & 0,0          & x+1, y-1\\
         x+1, y-1  & x-1, y+1 & 0,0
       \end{array}\right)
\end{equation}
Finally, let $\Gamma_{\alpha}$ denote the following game:
$$\Gamma_{\alpha}=\left(\begin{array}{cccc|c}
  C(x'_1,y'_1)  & 0                 & ... & 0               & [ 0 , y'_{1}-\alpha] \\
  0                 & C(x'_2,y'_2)  & ... & 0                   & [ 0 , y'_{2}-\alpha]\\
  ...               & ...               & ... & ...                       & ...         \\
 0              & 0                 & ... & C(x'_q,y'_q)    & [ 0 ,  y'_{q}-\alpha] \\
  \hline
\mbox{$[$}x'_{1}-\alpha , 0] & [x'_{2}-\alpha, 0] & ... & [x'_{q} -\alpha ,0]& G_{\alpha}
\end{array}\right)
$$
where $[x'_i-\alpha, 0]$ means a block of payoffs $(x'_i-\alpha,0)$ of appropriate size.

$\Gamma_0$ is a slight modification of the game obtained from $G_0$ by iterative applications ($q$ times) of lemma \ref{lm:add-cor}. 
Along the lines of the proof of lemma \ref{lm:add-cor}, it is easy to show that its Nash equilibria
correspond to the Nash equilibria of $G_0$ and that its extreme correlated equilibria are:
(i) its Nash equilibria;
(ii) the probability distributions with support in one of the blocks $C(x'_i,y'_i)$ and putting probability $1/6$
on every off-diagonal square of this block. It follows that $NEP(\Gamma_0)=U$ and $CEP(\Gamma_0)=P$.

Moreover, for any $\alpha$ small enough and any game $\Gamma$ in $B(\Gamma_{\alpha}, \alpha/2)$, every Nash equilibrium of $\Gamma_0$ is a Nash equilibrium of $\Gamma$; therefore, 
\begin{equation}
\label{eq:Gamma0-close-to-Gamma}
\forall (a_0, b_0)  \in NEP(\Gamma_0), \exists (a,b) \in NEP(\Gamma), ||(a, b)-(a_0, b_0)|| < 3\alpha/2
\end{equation}
(a closer look shows that we may replace $3\alpha/2$ by $\alpha/2$ in the RHS, but this is not needed).  Finally, let $\epsilon>0$. By upper semi-continuity of the Nash equilibrium correspondence, for $\alpha$ small enough, for any game $\Gamma$ in $B(\Gamma_0, 3\alpha/2)$,
\begin{equation}
\label{eq:Gamma-close-to-Gamma0}
\forall (a,b) \in NEP(\Gamma), \exists (a_0,b_0) \in NEP(\Gamma_0), ||(a_0,b_0)-(a,b)|| < \epsilon
\end{equation}
It follows from (\ref{eq:Gamma0-close-to-Gamma}) and (\ref{eq:Gamma-close-to-Gamma0}) that, for any $\alpha$ small enough and any game $\Gamma$ in $B(\Gamma_{\alpha}, \alpha/2)$, $NEP(\Gamma)$ and $NEP(\Gamma_0)$ are $\epsilon$-close. The same argument (up to replacement of Nash equilibrium by correlated equilibrium everywhere) shows that for every $\alpha$ small enough and for every game $\Gamma$ in $B(\Gamma_{\alpha}, \alpha/2)$, $CEP(\Gamma)$ and $CEP(\Gamma_0)$ are $\epsilon$-close.  Recalling that $NEP(\Gamma_0)=U$ and $CEP(\Gamma_0)=P$, this completes the proof.
\end{proof}

Note that, in the above proof, for every game $\Gamma$ in $B(\Gamma_{\alpha}, \alpha/2)$, any Nash equilibrium (resp. extreme correlated equilibrium) of $\Gamma_0$ is actually a strict Nash equilibrium (resp. strict correlated equilibrium\footnote{A correlated equilibrium $\mu$ is \emph{strict} if for every pure strategy $s_i$ with positive marginal probability under $\mu$, the inequalities in (\ref{eq:defcor}) are strict.}) of $\Gamma$. Since furthermore a convex combination of strict correlated equilibria is a strict correlated equilibrium, it follows that for every finite set $U$ in $\R^n$, every polytope $P$ containing $U$ and every $\epsilon>0$, there exists an open set of games $\Gamma$ such that: first, both the set of Nash equilibrium payoffs and the set of strict Nash equilibrium payoffs of $\Gamma$ are $\epsilon$-close to $U$; second, both the set of correlated equilibrium payoffs and the set of strict correlated equilibrium payoffs of $\Gamma$ are $\epsilon$-close to $P$.

\begin{appendix}
\section{Proofs of lemma \ref{lm:add-cor} and propositions \ref{prop:cor} and \ref{prop:gen}  in the $n$-player case.}

\textbf{Proof of lemma \ref{lm:add-cor}. } Let $G$ be a three-player game and let $(x,y,z) \in \R^3$. Assume w.l.o.g. that $x$, $y$, $z$, and all the payoffs in $G$ are positive. Let $C(x,y,z)$ denote the $3 \times 3 \times 1$ game (player $3$ is a dummy) where the payoffs of players $1$ and $2$ are as in (\ref{eq:Cxy}) and the payoff of player $3$ is always $z$. 
\begin{equation}
\label{eq:Cxyz}
C(x,y,z)=\left(\begin{array}{ccc}
         0,0, z          & x+1, y-1, z  & x-1, y+1,z \\
         x-1, y+1, z  & 0,0, z            & x+1, y-1, z \\
         x+1, y-1,z   & x-1, y+1, z   & 0,0, z
       \end{array}\right)
\end{equation}
Let
$$\Gamma= \hspace{0.1cm} \left(\begin{array}{c|c}
       C(x,y,z)
                            & [0,y,0]
        \\
        \hline
       [x,0,0] & [x,y,0]
     \end{array}\right)
\left(\begin{array}{c|c}
       [0,0,z] & [0,y,z]\\
        \hline
       [x,0,z] & G
     \end{array}\right)
$$
where $[x,0,0]$ denotes a block of payoffs $(x,0,0)$ of appropriate size. The same proof as in the two-player case shows that $\Gamma$ has the same set of NE payoffs as $G$ and that its set of CE payoffs is the convex hull of $(x,y,z)$ and of the set of CE payoffs of $G$. 

In the $n$-player case, the generalization of $C(x,y,z)$ simply consists in adding more dummy players with constant payoff. The generalization of $\Gamma$ should be clear from the following description of the three-player case: the pure strategy set of player $i$ is $S'_i \sqcup S_i$ (disjoint union), where $S_i$ and $S'_i$ are respectively the pure strategy set of player $i$ in $G$ and in $C(x,y,z)$. If for every $i$, player $i$ chooses a strategy in $S_i$ (resp. $S'_i$), then the payoffs are as in $G$ (resp. as in $C(x,y,z)$). If there exist players $i$ and $j$ such that $i$ chooses a strategy in $S_i$ and $j$ a strategy in $S'_j$, then the payoff of player $1$ (resp. $2$, $3$) is $x$ (resp. $y$, $z$) if he chooses a strategy in $S_1$ (resp. $S_2$, $S_3$) and $0$ otherwise.\\

\noindent \textbf{Proof of propositions \ref{prop:cor} and \ref{prop:gen}. } Let $U=\{ (x_i,y_i,z_i), 1 \leq i \leq m \}$ and let $P$ be the convex hull of $\{(x'_i,y'_i,z'_i), 1 \leq i \leq q \}$. Assume w.l.o.g. that the numbers $x_i$, $x'_i$, etc. are all positive. Let
$$G_{\alpha} = \left(
  \begin{array}{cccc|c}
    x_1,y_1,z_1  & 0 , 0, 0             &   ..              & 0, 0, 0                       & 0,y_1-\alpha,0\\
    0, 0, 0             & x_2, y_2, z_2 &   ..              & 0,0,0                             & 0,y_2-\alpha,0 \\
   ..                   & ..                    &  ..       & ..                                & ..            \\
    0,0,0              & 0,0,0                  & ..            & x_{m-1},y_{m-1} , z_{m-1}     & 0,y_{m-1}-\alpha,0\\
    \hline
x_1-\alpha,0,0 & x_2-\alpha,0,0 & ..            & x_{m-1}-\alpha, 0,0               & x_m , y_m, z_m \\
  \end{array}
\right)
$$
The same argument as in the two-player case shows that the Nash equilibria of $G_0$ are equal to its extreme correlated equilibria and correspond to the diagonal squares. It follows that $NEP(G_0)=U$ and that $CEP(G_0)$ is the convex hull of $U$. This proves proposition \ref{prop:cor} in the three-player case (for the $n$-player case, just add more dummy players). 

Now recall (\ref{eq:Cxyz}) and let $G'$ denote the $(3q+m) \times (3q+m)\times 1$ game with block diagonal payoff matrix
$$G'=\left(
  \begin{array}{ccc}
    C(x'_1,y'_1,z'_1)   &  & 0  \\
                      & \ddots &        \\
    0                       &  & C(x'_q,y'_q,z'_q)   \\
    \end{array}
\right) $$
Let $x$, $y$, $z$ be positive real numbers. Let $\Gamma_{\alpha}$ denote the $(3q+m)\times (3q+m)\times 2$ game:
$$\Gamma_{\alpha}=
\left( \begin{array}{c|c}
            G'       &  [0, y'_i-\alpha, 0]  \\
\hline
    [x'_j-\alpha, 0, 0]   &  [x, y,0] \\
  \end{array}\right)
  \left(\begin{array}{c|c}
  [0, 0, z'_i-\alpha]    & [0,y,z] \\
\hline
[x,0, z] & G_{\alpha}\\
\end{array}\right)$$
This should be read as follows: if player 1 chooses row $i>3q$ and player $2$ chooses column $j \leq 3q$ with $3p+1\leq j \leq 3p+3$, then the payoffs are $(x'_p-\alpha,0,0)$ if player $3$ chooses the left matrix and $(x,0,z)$ if player $3$ chooses the matrix on the right.

Fix $\epsilon>0$ and assume that $P \supset U$. The same arguments as in the two-player case show that, for every $\alpha$ small enough, and for every game $\Gamma$ in $B(G_{\alpha}, \alpha/2)$, $NEP(\Gamma)$
is $\epsilon$-close to $U$ and $CEP(\Gamma)$ is $\epsilon$-close to $P$. The $n$-player case is similar.\footnote{Instead of $C(x,y,z)$, $\Gamma$, $G_{\alpha}$, $\Gamma_{\alpha}$, and their $n$-player version, it is possible to use games in which the roles of the players are symmetric, but this is less parcimonious.} 
\end{appendix}


\begin{thebibliography}{1}
\bibitem{Aumann} Aumann, R. (1974), ``Subjectivity and Correlation in Randomized Strategies'',
\emph{Journal of Mathematical Economics} \textbf{1}, 67-96
\end{thebibliography}
\end{document}